\documentclass{article}

\usepackage{arxiv}

\usepackage[utf8]{inputenc} 
\usepackage[T1]{fontenc}    
\usepackage{hyperref}       
\usepackage{url}            
\usepackage{booktabs}       
\usepackage{amsfonts}       
\usepackage{nicefrac}       
\usepackage{microtype}      
\usepackage{lipsum}
\usepackage{graphicx}
\graphicspath{ {./images/} }

\title{Product-Form Distribution and Reversibility of Inhomogeneous Symmetric Simple Exclusion Process with Open Boundaries
}

\author{
Marina V. Yashina  \\
  Department of Higher Mathematics\\
  Moscow Automobile and Road Construction\\
  State Technical University (MADI) \\
  Moscow, Leningradsky avenue, 64, Russia  \\
  \texttt{mv.yashina@madi.ru} \\ 
\And
 Alexander G. Tatashev  \\
  Department of Higher Mathematics\\
  Moscow Automobile and Road Construction\\
  State Technical University (MADI) \\
  Moscow, Leningradsky avenue, 64, Russia  \\
  \texttt{a-tatashev@yandex.ru} \\
}

\begin{document}
\maketitle
\begin{abstract}
 
 We consider an inhomogeneous symmetric simple exclusion process on a one-dimensional lattice with open boundary conditions. The time scale is continuous. Particles of different types arrive to the utmost left and the utmost right site. If a particle is in a site that is neither the utmost left site nor the utmost right site, then the particle moves onto one site to the left or to the right. If a particle is either in the utmost left site  or the utmost right site, then the particle leaves the system or moves onto one cell to the right or to the left, respectively. An arrival or a transfer of particle is possible only to a vacant site. The rate of arrival, exit or movement of a particle depends on its type and does not depend on the site from that the particle arrives or exits and on the direction the movement. The stationary distribution of the system states probabilities has been found. This distribution turns out to be multiplicative in the sense that the probability of the site state does not depend on the states of the other sites in the stationary mode, and the steady probability of any state of the system is equal to the product of the site states steady probabilities, and the probability of site ocupancy is the same as the server occupancy probability for M/G/1/0 loss system. We have found the arrival rate and the average sojourn time have been found.  We have proved the reversibility of the process in time under the additional condition that the rate of arrival of a particle of prescribed type is equal to the rate of this type particle departure.

\end{abstract}


\section{Introduction}
An exclusion process is a model such that particles move along a lattice divided into sites. No more than one particle can be in the same site simultaneously. A particle can only move  to a vacant site. Exclusion processes were first considered in [1], where exclusion processes with continuous time were studied. In the case of the continuous time exclusion process, a particle atempts to move after a time interval with exponentially distributed duration, and the attempt to move is realized under the condition that the site to that the particle atempts to move is vacant. For exclusion processes with a discrete time scale, it is possible for more than one particle to move at the same time (a process with parallel updating or a synchronous exclusion process).  This complicates the study of such processes. Discrete-time exclusion processes with sequential updating (asynchronous exclusion processes) are also considered. For these processes only one particle can move at any discrete moment. The results obtained in the study of discrete-time exclusion processes with sequential updating are analogous to the results obtained in the study of continuous-time exclusion processes, and therefore continuous-time exclusion processes are also called exclusion processes with sequential updating or asynchronous exclusion processes, and usually only continuous-time asynchronous exclusion processes are considered, since the results are easily transferred to the related asynchronous discrete-time processes. An  exclusion process on a one-dimensional lattice is called simple if a particle can move only to a nearest neighboring site. Simple exclusion processes on a finite closed lattice (lattice with periodic boundary conditions), a finite lattice with open boundary conditions, or an infinite lattice were considered. In [1],  zero-range processes were also introduced into consideration. These processes differ from exclusion processes in that any number of particles can be in the same site. Models of statistical physics and models of traffic flows are based on exclusion processes and zero-range processes. Exclusion processes were considered, in particular, in [1]--[35]. In [31], a continuous-time exclusion process  was considered such that there is a finite set of particles on an infinite one-dimensional lattice, and any particle hops onto one site to the right or to the left under the assumption that the related site is vacant, and the hop attempt rate depends on the particle.

Exclusion processes and zero-range processes of some types are equivalent to processes related to queueing networks. Known results of the queueing theory useful for exclusion processes study are provided, e.~g., in [5], [9], [14], [29], [31].  One of subject of queueing network study is to find conditions of multiplicativity, i.e., the possibility of representing probabilities of network states as a product of stationary probabilities of node states (this is equivalent to the fact that the probability of a node state does not depend on the probabilities of the states of the remaining nodes).  Another subject is to find conditions of reversibility, i.~e., the equivalence of the original process and a process that differs from the original in the direction of the time axis, [36]--[38]. The same questions were considered for exclusion processes and zero-range processes,  in particular, in [1], [4]--[6], [10], [14], [28], [29], [31]. 

In this paper, we consider a continuous-time exclusion process on a lattice containing a finite number of sites. Particles of different types enter and leave the system through the utmost left and utmost right sites. A particle can enter a site only if the site is vacant at the moment of the particle arrival. The rates of particle arrival and depature depend on the particle type but do not depend on whether the arrival or exit occurs at the leftmost or rightmost site. A particle can move to adjacent site on the left or on the right under the condition that this site is vacant, and the rate of movement depends on the particle type but not on the direction of movement. In this sense, the process under consideration is symmetric. In Section 2, we describe the process. In Section~3, we prove the ergodicity of the process. Namely, it is proved that there exist the steady probabilities of process states that are independent of the initial state. In Section~4, we have found steady probabilities of system states and the steady probabilities of site occupancy. It is proved that the steady probability of system state is represented as the product of the steady probabilities of sites (product-form distribution), and the steady probability of site ocupancy is the same as the server occupancy probability for M/G/1/0 loss system. In Section 5, we have ontained formulas for computing the arrival rates and the average sojourn time of a prescribed particle. In Section 6, we consider the process  under the additional condition that, for each type, the rate of arrival of prescribe type particles is equal to the rate of departure of this type particle. We have shown that under this assumption, all states of the process are equiprobable in the stationary mode, and the process is reversible.

\section{
Description of system
}
\label{section:Description} Suppose the system contains a one-dimensional lattice divided into $N\ge 2$~sites. At any time $t\in [0,+\infty),$ in a site, there is  a particle of one of $K$ types or the site is vacant. The sites are numbered from the left to the right. The indices of the sites are $i=1,\dots,N.$ If, at time $t,$ the site~1 (leftmost) or site~$N$ (rightmost) is vacant, then, in the time interval $(t,t+\Delta t),$ with probability $\alpha_k\Delta t+o(\Delta t),$ $\Delta t\to 0,$ a particle of type~$k$ enters this site from outside and, at time $t+\Delta t,$ this particle will be in this site. If at time $t$, in site~1 or in site~$N,$ there is a particle of type~$k,$ then, in the time interval $(t,t+\Delta t),$ with probability
$\beta_k\Delta t+o(\Delta t),$ $\Delta t\to 0,$ the particle will leave the system, and, at time $t+\Delta t,$ this site will be vacant. If, at time~$t,$ a particle of the type $k$ is in the site~$i,$ and the site $i-1$ (the site $i+1)$ is vacant  then, in the time interval $(t,t+\Delta t),$ with probability $\delta_k\Delta t+o(\Delta),$ $\Delta t\to 0,$ the particle hops to the site $i-1$ (site $i+1),$ $i=2,\dots,N-1.$

\section{
System ergodicity
}
\label{section:SystEr} Suppose $x=(x_1,\dots,x_N)$ is the state of the system such that $x_i$ is the state of the site~$i;$ $x_i=0$ if the site~$i$ is vacant and $x_i=k$ if the site~$i$ is occupied by a particle of the type~$k,$ $k=1,\dots,K;$ $x(t)=(x_1(t),\dots,x_N(t))$ is the stochastic process, where $x(t)$ is the system state at  time~$t,$ $k=1,\dots,K,$  $i=1,\dots,N.$ The stochastic process $x(t)$ is a continuous-time Markov chain. Denote by $G$ the process $x(t)$ state space. The space $G$ consists of $(K+1)^N$ states $(x_1,\dots,x_N)$ such that $0\le x_i\le K$ are intergers, $i=1,\dots,N.$ 
\vskip 3pt
{\bf Theorem 1.} {\it The Markov chain $x(t)$ is ergodic. Namely, there exists the steady probability  of any state $(x_1,\dots,x_N)\in G,$ which does not depend on the initial state, and this probability is not equal to~zero. 
}
\vskip 3pt
According to an ergodic theorem for a continuos-time Markov chain [39], if the state space is finite, and the space consists of a unique communicating state class (i.~e., the Markov chain is irreducible), then there exists a non-zero  probability of any state independent of the initial state.  The Markov chain~$x(t)$ satisfies this condition.  Indeed, the system can transfer from any state to any another state. Indeed, from any state $x\in G$, in time with finite expectation, the process $x(t)$ can transfer to the state $(0,\dots,0),$ and from the state $(0,\dots,0),$ in time with finite expectation, the process can transfer to any other state. Theorem~1 has been proved.
\vskip 3pt
Suppose $p(x_1,\dots,x_N)$ is the steady probability of the state $(x_1,\dots,x_N),$ and $P_i(x_i)$ is the steady probability that the site $i$ is in the state $x_i,$ $i=1,\dots,N.$

\section{
Steady state probabilities
}
\label{section:SteadySt}

In Section 4, we prove that the stationary state distribution of the process $x(t)$ is a product-form distribution. 
\vskip 3pt
{\bf Theorem 2.} {\it Steady probabilities of site states are computed as 
$$p(x_1,\dots,x_N)=C\prod\left(\frac{\alpha_{x_i}}{\beta_{x_i}}\right),\eqno(1)$$ 
where the product is taken over $i$ such that $x_i\ne 0;$
$$C=p(0,\dots,0)=\frac{1}{\left(1+\sum_{k=1}^K\limits\frac{\alpha_k}{\beta_k}\right)^N}.\eqno(2)$$
}
\vskip 3pt 
{\bf Proof.} Denote by $M$ the number of the process $x(t)$ states $(M=(K+1)^N).$ Suppose the process $x(t)$ states are numbered arbitrary; $p_i$ is the steady state probability of the state~$i;$ $\lambda_{ij}$ is the rate of hop from the state~$i$ to the state~$j;$ $i,j=1,\dots,M,$ where $M$ is the number of the process $x(t)$ states $(M=(K+1)^N).$ The system of equations for steady state probabilities
$$p_i\sum\limits_{j=1}^{M}\lambda_{ij}=\sum\limits_{j=1}^{M}p_j\lambda_{ji},\ i=1,\dots,M,\ \eqno(3)$$
$$\sum\limits_{i=1}^{M}p_i=1.\eqno(4)$$

Let $i$ and $j$ be states such that $\lambda_{ij}>0.$ The transition from the state $i$ to the state $j$ occurs if a particle enters the system, or a particle leaves the system, or a particle passes from one state to another state. 

Suppose $A_k$ is the set of pairs $(i,j)$ such that the transition from the state~$i$ to the state~$j$ is due that  a particle of the type~$k$ particle arrives into the system;  $B_k$ is the set of pairs $(i,j)$ such that the transition from the state~$i$ to the state~$j$ occurs if a particle of type~$k$ leaves the system; $D_k$ is the set of pairs $(i,j)$ such that the transition from the state~$i$ to the state~$j$ occurs if  a particle of the type~$k$ moves. 

Suppose $(i,j)$ belongs to $A_k.$ Then $(j,i)$ belongs to $B_k,$
$$\lambda_{ij}=\alpha_k,\ (i,j)\in A_k,\ k=1,\dots,K,\eqno(5)$$
$$\lambda_{ji}=\beta_k,\ (i,j)\in A_k,\ k=1,\dots,K.\eqno(6)$$
In this case, in accordance with (1), we have
$$\frac{p_j}{p_i}=\frac{\alpha_k}{\beta_k},\ (i,j)\in A_k,\ k=1,\dots,K.\eqno(7)$$
From (5)--(7), we obtain
$$p_i\lambda_{ij}=p_j\lambda_{ji},\ (i,j)\in A_k,\ k=1,\dots,K.\eqno(8)$$
From (8), it follows that any term on the left-hand side of (3) such that $(i,j)\in A_k$ corresponds to an equal term on the right-hand side of (3).

Suppose that $(i,j)\in B_k.$ Then $(j,i)\in A_k;$
$$\lambda_{ij}=\beta_k,\ (i,j)\in B_k,\ k=1,\dots,K,\eqno(9)$$
$$\lambda_{ji}=\alpha_k,\ (i,j)\in B_k,\ k=1,\dots,K,\eqno(10)$$
$$\frac{p_j}{p_i}=\frac{\beta_k}{\alpha_k},\ (i,j)\in B_k,\ k=1,\dots,K.\eqno(11)$$
From (9)--(11) it follows that
$$p_i\lambda_{ij}=p_j\lambda_{ji},\ (i,j)\in B_k,\ k=1,\dots,K.\eqno(12)$$
As follows from (12), any term in the left-hand side of (3) such that $(i,j)\in B_k$ corresponds to a term with the same value in the right-hand side~(3).

Let $(i,j)\in D_k.$ Then $(j,i)\in D_k;$
$$\lambda_{ij}=\delta_k,\ (i,j)\in D_k,\ k=1,\dots,K,\eqno(13)$$
$$\lambda_{ji}=\delta_k,\ (i,j)\in D_k,\ k=1,\dots,K,\eqno(14)$$
and, in accordance with (14),
$$p_i=p_j,\ (i,j)\in D_k,\ k=1,\dots,K.\eqno(15)$$
From (13)--(15) it follows that
$$p_i\lambda_{ij}=p_j\lambda_{ji},\ (i,j)\in D_k,\ k=1,\dots,K.\eqno(16)$$
From (16) it follows that any term on the left-hand side of (3) such that $(i,j)\in D_k$ corresponds to a term on the right-hand side of~(3) with the same value.

Thus any term on the left-hand side of~(3) corresponds to an
equal term on the right-hand side of~(3). This is a one-to-one correspondence.
Indeed, a left-hand side term related to the pair $(i,j)$ corresponds to a right-hand side term related to the pair $(j,i)$ and thus different terms on the left-hand side of~(3) correspond to different terms on the right-hand side of (3), and the number of non-zero terms in (3) is the same on the left-hand and right-hand sides, since $\lambda_{ij}>0$ if and only if $\lambda_{ji}>0.$ From the established one-to-one correspondence of the terms on the left-hand and the right-hand sides of (3,) it follows that (1), (2) satisfies (3).

Substituting (1) into the left-hand side of (4), we can represent the equation for $C$ in the form
$$C\left(1+\sum_{k=1}^K\limits\frac{\alpha_k}{\beta_k}\right)^N=1.$$
Thus (1), (2) is a solution to system (3), (4).
Theorem~2 has been proved.
\vskip 3pt

{\bf Theorem 3.} {\it The stationary state distribution is the same  for any site, and the state probabilities are computed as
$$P_i(0)=\frac{1}{1+\sum\limits_{k=1}^K\frac{\alpha_k}{\beta_k}},\ i=1,\dots,N, 
\eqno(17),$$
$$P_i(k)=\frac{\alpha_k}{\beta_k}\left(1+\sum\limits_{l=1}^K\frac{\alpha_l}{\beta_l}\right)^{-1},\ k=1,\dots,K,\ i=1,\dots,N.\eqno(18)$$
}
\vskip 3pt 
Theorem 3 follows from Theorem 2. 
\vskip 3pt 
Note that, in accordance  the Erlang B formula, the value computed according to~(18) is the steady state that a type $k$ request is in service in a one-channel loss system with an inhomogeneous Poisson input, where $\alpha_k$ is the type~$k$ request arrival rate, and $1/\beta_k$ is the average service time. 
\vskip 3pt 
{\bf Theorem 4.} {\it The stochastic process $x(t)$ has a product-form distribution 
$$p(x_1,\dots,x_N)=\prod_{i=1}^NP_i(x_i),\ (x_1,\dots,x_N)\in G.$$ 
}
Theorem 4 follows from Theorems 2, 3.

\section{Arrival rates and average sojourn time for prescribed type particle}
\label{section:RelSt}

The average number of the type $k$ particles arriving to the system per a time unit in stationary mode is called the type~$k$ particle arrival rate. Denote by $J_k$ this rate,  $k=1,\dots,k$ 
\vskip 3pt
{\bf Theorem 5.}  {\it The type~$k$ particle arrival rate is computed as
$$J_k=\frac{2\alpha_k}{1+\sum\limits_{l=1}^K\frac{\alpha_l}{\beta_l}},\ k=1,\dots,K.\eqno(19)$$
}
\vskip 3pt
{\bf Proof.} We have
$$J=\alpha_k (P_1(0)+P_N(0)).\eqno(20).$$
Combining (17) and (20), we obtain (19). Theorem 5 has been proved.
\vskip 3pt
Denote by $U_k$ the average sojourn time for a particle of type~$k.$
\vskip 3pt
{\bf Theorem 6.} {\it The following formula holds 
$$U_k=\frac{N}{2\beta_k},\ k=1,\dots,K.\eqno(21)$$
}
\vskip 3pt
{\bf Proof.} According to Little's formula, the average sojourn time  for a prescribe type request in a queueing system equals the ratio of the average number of requests of this type in the system to the arrival rate for the requests of this type arrival. The system under consideration can be interpreted as a queueing system, where requests correspond to particles of related types. The average number of particles of type~$k$ in the system is $P_1(k)+\dots+P_N(k).$ Taking into account (18), (19) and Little's formula, we obtain (21). Theorem~6 has been proved.
\vskip 3pt
The average number of the type $k$ particles arriving to the system per a time unit in stationary mode is called the type~$k$ particle arrival rate. Denote by $J_k$ this rate,  $k=1,\dots,k$ 
\vskip 3pt
{\bf Theorem 5.}  {\it The type~$k$ particle arrival rate is computed as
$$J_k=\frac{2\alpha_k}{1+\sum\limits_{l=1}^K\frac{\alpha_l}{\beta_l}},\ k=1,\dots,K.\eqno(19)$$
}
\vskip 3pt
{\bf Proof.} We have
$$J=\alpha_k (P_1(0)+P_N(0)).\eqno(20).$$
Combining (17) and (20), we obtain (19). Theorem 5 has been proved.
\vskip 3pt
Denote by $U_k$ the average sojourn time for a particle of type~$k.$
\vskip 3pt
{\bf Theorem 6.} {\it The following formula holds 
$$U_k=\frac{N}{2\beta_k},\ k=1,\dots,K.\eqno(21)$$
}
\vskip 3pt
{\bf Proof.} According to Little's formula, the average sojourn time  for a prescribe type request in a queueing system equals the ratio of the average number of requests of this type in the system to the arrival rate for the requests of this type arrival. The system under consideration can be interpreted as a queueing system, where requests correspond to particles of related types. The average number of particles of type~$k$ in the system is $P_1(k)+\dots+P_N(k).$ Taking into account (18), (19) and Little's formula, we obtain (21). Theorem~6 has been proved.

\section{Equiprobability of states and reversibility under the condition
$\alpha_k=\beta_k$}
\label{section:Equip}

In Section 5, we consider the system under an assumption that
$$\alpha_k=\beta_k,\  k=1,\dots,K.\eqno(22)$$ 
\vskip 3pt
{\bf Theorem 7.} {\it If (31) holds, then the steady state is the same for any state. Namely,
$$p(x_1,\dots,x_N)=\frac{1}{(K+1)^N},\eqno(23)$$ 
$x_1,\dots,x_N$ are intergers,  $0\le x_1,\dots,x_N\le K.$
}
\vskip 3pt
{\it Rroof.} Combining (1), (2) and (22) , we obtain (23). Theorem~7 has been proved.
\vskip 3pt 
Further we will consider the process $x(t)$ on the infinite in both directions time interval $(-\infty,\infty).$ We will also consider on the time interval $(-\infty,\infty)$ a random process $\tilde{x}(t)=x(-t)$ that differs from the process $x(t)$ in the direction of the time axis. Suppose $P(x(t_1)=E_j|x(t_0)=E_i)$ is the probability that, at time $t_1,$ the system will be in the state $j$ under the assumption that, at time $t_0,$ the system was in the state~$i;$ $\tilde{\lambda}_{ij}$ is the rate of the transition of the process $\tilde{x}(t)$ from the state~$i$ to the state~$j.$
\vskip 3pt
{\bf Theorem 8.} {\it Stochastic processes $x(t)$ and $x(-t)$ are stochastic equivalent.}
\vskip 3pt
{\bf Proof.} Suppose $(i,j)\in A_k$ or $(i,j)\in B_k$
We have, $\Delta t\to 0,$
$$P(x(-t-\Delta t)=E_j|P(x(-t)=E_i)=\frac{P(x(-t-\Delta t)=E_j)\alpha_k\Delta t+o(\Delta t)}{P(x(-t)=E_i)}.\eqno(24)$$
Taking into account Theorem 1, (12), and (24), we have
$$\tilde{\lambda}_{ij}=\alpha_k=\lambda_{ji}=\lambda_{ij}\eqno(25)$$
under the assumption that the transition from the state $j$ to the state $i$ 
(and therefore the transition of the process $\tilde{x}(t)$ from the state~$i$ to the  state $j)$ due to an arrival or a departure of the type $k$ particle.

Similary, we have
$$P(x(-t-\Delta t)=E_j|P(x(-t)=E_i)=\frac{P(x(-t-\Delta t)=E_j)\delta_k\Delta t+o(\Delta t)}{P(x(-t)=E_i)}$$
$$\tilde{\lambda}_{ij}=\delta_k=\lambda_{ij}\eqno(26)$$
if the transition of the process $x(t)$ from the state $j$ to the state $i$ (and, consequently, the transition of the process $x(t)$ from the state~$i$ to the state $j)$ is due to a hop of a particle of type~$k;$
$$P(x(-t-\Delta t)=E_j|P(x(-t)=E_i)=o(\Delta t),$$
$$\tilde{\lambda}_{ij}=0=\lambda_{ij}\eqno(27)$$
if the transition in one step from the state $i$ to the state $j$ is impossible.

Theorem 8  follows from (25)--(27).

\section{Conclusion}
\label{section:Concl}

We have proved that an inhomogeneous symmetric continuous one-dimensional exclusion process with open boundary conditions has a product-form stationary distribution under the assumption that the process parameters are symmetric, and the formula for steady probabilities has a simple form. Namely, it is proved that the steady probability that a site is occupied equals the steady probability that there is a request in service in M/G/1/0 loss system.  Knowledge of the stationary distribution allows us to compute the steady probability of site occupancy, the arrival rate, and the average sojourn time for a particle of prescribed type.  Under an additional assumption, the process reversibility has been proved.

\bibliographystyle{unsrt}  

{}

\end{document}